\documentclass[12pt, oneside]{article}   	
\usepackage{geometry}                		
\usepackage{graphicx}				
\usepackage{amssymb}

\def\E{{\rm I\!E}}
\def\R{{\rm I\!R}}
\def\pos{{\rm pos}}
\def\c.#1{{\cal #1}}
\def\b.#1{{\bf #1}}
\def\QED{\hfill$\bullet$}

\title{\bf ?`How do 9 points looks like in $\E^3$?}
\author{\sc Ricardo Strausz  \\ Instituto de Matem\'aticas, U.N.A.M. \\ \tt dino@math.unam.mx}
\date{}							

\begin{document}
\maketitle
\abstract{The aim of this note is to give an elementary proof of the following fact: given 3 red convex sets and 3 blue convex sets in $\E^3$, such that every red intersects every blue, there is a line transversal to the reds or there is a line transversal to the blues. This is a special case of a theorem of Montajano and Karasev \cite{MK} and generalizes, in a sense, the colourful Helly theorem due to Lov\'asz (cf. \cite{BL})}

\section{Introduction.}

Consider 9 points in the euclidian 3-dimensional space $\E^3$, and order them, by free will, in a $3\times3$ matrix:

$$\bordermatrix{
\E^3 & \c.U & \c.V & \c.W \cr
\c.A & a & \alpha & x \cr
\c.B & b & \beta & y \cr
\c.C & c & \gamma & z \cr
}$$

Now, consider the convex hull of the triples of its rows, say 3 blue triangles $\c.A, \c.B$ and $\c.C$, and of its columns, say 3 red triangles $\c.U, \c.V$ and $\c.W$; observe that each red triangle intersect each blue triangle.

The aim of this note is to exhibit an elementary proof of the following fact:

\medskip
\noindent
{\bf Theorem.} {\it There is a line transversal to the red triangles, or there is a line transversal to the blue triangles.\/}
\medskip

This was first observed by Luis Montejano \cite{M} and later generalised in collaboration with Roman Karasev \cite{MK}. Their proof uses multiplication formulas for Schubert cocycles, the Lusternik-Schnirelmann category of the Grassmannian, different versions of the colorful Helly theorem by B\'ar\'any and Lov\'asz, and a bit of Separoid Theory (see \cite{ABMS,S}). In contrast, here I only use very elementary facts to prove this case: the non-planarity of the graph $K_{3,3}$, the fact that given $d+2$ convex sets in $\R^n$ either they admit a $d$-transversal or every subset of them is separated from its complement, and the following result ---which is interesting on its own--- which can be proved using only basic properties of the interior product (Section~3).

\medskip
\noindent
{\bf Basic Lemma.} {\it Let $H_A^+, H_U^+, H_W^+$ and $H_C^+$ be four affine semispaces in $\E^n$, with unitary normal vectors $A,U,W$ and $C$, respectively. If
	$$(H_A^+\cap H_U^+)\setminus(H_W^+\cup H_C^+)\not=\emptyset,$$
and
	$$(H_W^+\cap H_C^+)\setminus(H_A^+\cup H_U^+)\not=\emptyset,$$
then
	$$\pos(A,U)\cap\pos(W,C)=\emptyset,$$
\/}
Here, and in the sequel, $\pos()$ denotes the positive span.

\section{Proof of Theorem.}

Consider 3 red triangles and 3 blue triangles in $\E^3$ as in the Introduction (each red intersects each blue), and suppose that there are no line transversals neither to the red, nor the blue. Then, each red triangle is separated by a plane from the other two red triangles, and each blue triangle is separated too --- each triangle must be disjoint from the convex hull of the other two triangles of the same color. Let us denote by $H_A^+$ the semispace, with (blue) normal vector $A$, that is the witness of the separation $\c.A\mid\c.B\c.C$; that is, $\c.A\subset H_A^+$ and $(\c.B\cup\c.C)\cap H_A^+=\emptyset$. Analogously, we consider the witnesses of all separations, and take the 3 blue vectors $A,B$ and $C$, and the 3 red vectors $U,V$ and $W$ in the unitary sphere of $\E^3$. 

Now, join each blue vector, with a spherical segment, to each red vector and, abusing the notation, name the segment with that point of the original configuration that corresponds to the intersection of the respective triangles; e.g., the unitary vectors $A$ and $U$ are joined with the spherical segment $a$ while the unitary vectors $W$ and $C$ are joined with the spherical segment $z$. We have just drawn $K_{3,3}$ in the 2-sphere, therefore there must be a crossing, which we may suppose is the crossing of the spherical segments $a$ and $z$. 

However, we have that our original points $a$ and $z$ satisfy
	$$a\in(H_A^+\cap H_U^+)\setminus(H_W^+\cup H_C^+),$$
and
	$$z\in(H_W^+\cap H_C^+)\setminus(H_A^+\cup H_U^+),$$
thus the crossing contradicts the Basic Lemma. \QED

\section{Proof of Lemma.}

Let 
	$$a\in(H_A^+\cap H_U^+)\setminus(H_W^+\cup H_C^+),$$ 
	$$z\in(H_W^+\cap H_C^+)\setminus(H_A^+\cup H_U^+),$$
and suppose 
	$$\eta\in\pos(A,U)\cap\pos(W,C).$$
\noindent
Then, there exist $i,j,r,s>0$ such that $\eta=iA+jU=rW+sC$.

We may suppose that $0\in H_A\cap H_U$ and $p\in H_W\cap H_C$ ---here $H_A$ denotes the hyperplane bounding $H_A+$, and respectively with $H_U$ and the others. Then, the hypothesis of the lemma can be rewritten as the following eight inequalities:
$$\matrix{
a\cdot A>0 & z\cdot A<0 \cr
a\cdot U>0 & z\cdot U<0 \cr
(a-p)\cdot W<0 & (z-p)\cdot W>0 \cr
(a-p)\cdot C<0 & (z-p)\cdot C>0\cr
}$$

But, on the one hand we have that $\eta\cdot a=(iA+jU)\cdot a>0$ and $\eta\cdot(a-p)=(rW+sC)\cdot(a-p)<0$, therefore
	$$0<\eta\cdot a<\eta\cdot p;$$
on the other hand we have that $\eta\cdot z=(iA+jU)\cdot p<0$ and $\eta\cdot(z-p)=(rW+sC)\cdot(z-p)>0$, therefore
	$$0>\eta\cdot z>\eta\cdot p,$$
a clear contradiction which arises from the supposition of the existence of $\eta$.\QED

\section{Remarks and Open Problems.}

It is easy to see that the basic lemma can be generalised to arbitrary large hyperplane pencils ---the proof of this is totally analogous to that in Section~3.

\bigskip
\noindent
{\bf Lemma\/} {\it Consider two pencils of semispaces $\c.R=\{H^+_{v_i} : i=1,\dots,r\}$ and $\c.S=\{H^+_{u_j} : j=1,\dots,s\}$, with normal unitary vectors $\{v_1,\dots,v_r\}$ and $\{u_1,\dots,u_s\}$, respectively. If 
$$\cap\c.R\setminus\cup\c.S\not=\emptyset\ {\rm and }\ \cap\c.S\setminus\cup\c.R\not=\emptyset,$$ 
then 
$$\pos(v_1,\dots,v_r)\cap\pos(u_1,\dots,u_s)=\emptyset.$$\QED
}

\bigskip
But, how can we use this to prove the general case? What about the special case $4\times3$

$$\bordermatrix{
\E^4 & \c.U & \c.V & \c.W \cr
\c.A & a & \alpha & x \cr
\c.B & b & \beta & y \cr
\c.C & c & \gamma & z \cr
\c.D & d & \delta & w \cr
}?$$

\noindent
Here the conclusion is that: {\sl there exists a 2-plane transversal to the four triangles, or a line transversal to the three tetrahedra\/} ---recall that, given 4 convex sets, either there is a 2-plane transversal to those 4 sets, or each of them is separated from the other three and each pair is separated from the other two (cf. the definition of {\sl Radon dimension\/} in \cite{S}). 

For, consider the uniform rank 3 hypergraph $K^3_{4,3}=(R\cup B,\Delta)$ consisting of seven vertices, four blue and three red, and all those triangles which can be formed by taking two blue vertices and one red. Observe that $K^3_{4,3}$ cannot be embedded into the 3-sphere without crossings; the link of each blue vertex is a drawing of a $K_{3,3}$ in the 2-sphere.

\medskip
There are other two lines of research which deserves attention: 

---Is the Theorem true in the context of oriented matroids? That is, can it be settled for all oriented matroids of order 9 and rank 4? (In this context, we may need to change the notion of {\sl line transversal to the coloured triangles\/} by {\sl the convex hull of a coloured triangle is not separated from the convex hull of the others\/}.)

---Is the basic lemma true in such a context? That is, can we change ``semispaces'' by ``pseudosemispaces''?

The first question seems easy to check ---there are 9,276,601 uniform oriented matroids to be checked--- while the second may deserve a deeper understanding of the topological representation theorem.

\end{document}